\documentclass[11pt]{amsart}
\numberwithin{equation}{section}
\newtheorem{theorem}{Theorem}

\numberwithin{theorem}{section} \numberwithin{lemma}{section}
\numberwithin{proposition}{section}

\def\al{\aligned}
\def\eal{\endaligned}

\begin{document}

\tracingpages 1
\title[gradient estimate]{\bf Some gradient estimates for the heat
equation on domains and for an equation by Perelman}
\author{ Qi S. Zhang}
\address{Department of Mathematics,  University of California,
Riverside, CA 92521, USA }
\date{May 2006}

\begin{abstract}
 In the first part,  we  derive a sharp gradient estimate for the log of Dirichlet
  heat kernel and Poisson heat kernel on domains, and a sharpened local
  Li-Yau gradient estimate.

  In the second part, without explicit curvature
  assumptions, we prove a global upper bound for the fundamental
  solution of an equation introduced by G. Perelman, i.e. the heat
  equation of the conformal Laplacian under backward Ricci flow.
  Further, under nonnegative Ricci curvature assumption,
  we prove a qualitatively sharp, global Gaussian upper
  bound. The idea is to combine the Nash and Davies heat kernel estimate
  with a Sobolev imbedding by Hebey, together with a Hamilton type
  gradient estimate.

\end{abstract}
\maketitle
\tableofcontents
\section{Introduction}

The goal of the paper is to establish certain new point-wise or
gradient estimates for the heat equation in both the fixed metric
and the Ricci flow case. Gradient estimates for the heat equation
are important tools in geometric analysis as pointed out in the
papers [LY], [H], [CH], [ACDH] and others. In this paper, for the
fixed metric case, two gradient estimates are proven. One is a sharp
gradient estimate for the log of Dirichlet heat kernel and Poisson
heat kernel on domains. This can be viewed as a boundary version of
the well known Li-Yau gradient estimate (see Section 2 for a
restatement).   It can also be viewed as another step in the long
running process of heat kernel or Poisson kernel estimate starting
with the Gaussian formula and Poisson formula. As far as boundary
gradient estimate is concerned, only the Neumann boundary case was
treated in [LY] and [Wa]. The Dirichlet case is different in that
solutions vanish on the boundary. Therefore, the estimate is
different from both the Li-Yau theorem and its generalization in
[Wa]. The result seems to be new even for Euclidean domains. The
other result is a sharpened local Li-Yau gradient estimate that
matches the global one. As well known, even for manifolds with
nonnegative Ricci curvature, the local Li-Yau estimate differs with
the global one by a nontrivial factor. Here we show that this factor
can be chosen as one. Hence the global and local Li-Yau estimate are
identical. We expect the result to have applications in heat kernel
estimate on manifolds. These two results are presented in section 2.

In sections 3 and later we turn to the case when the metric
evolves by the Ricci flow. In the fundamental paper [P], Perelman
discovered a monotonicity formula for equation (4.0) below, which
can be regarded as the heat equation of the conformal Laplacian
under backward Ricci flow. Perelman's formula can be thought of as
a gradient estimate. Using this estimate together with the reduced
distance function, he then applies the maximum principle to prove
a lower bound for the fundamental solution of (4.0), whenever it
exists. The reduced distance incorporates the scalar curvature as
an integral part. However it does not seem that the maximum
principle alone will yield a two sided bound.  Here we will use
the Nash method to show that certain long time upper bound also
holds. The bounds involve more classical geometric quantities such
as the best constants in Sobolev imbedding, which depend only on
Ricci curvature lower bound and injectivity radius. Under no
explicit curvature assumption, we prove a global on-diagonal upper
bound for the fundamental solution on (4.0). The bound is good in
the sense that it matches the on diagonal bound in the fixed
metric case. However, we are not able to prove a good off-diagonal
bound without further assumptions on curvature. Nevertheless, the
on-diagonal upper bound does not have the usual, trouble making
exponentially growing term even when the Ricci curvature changes
sign. When the Ricci curvature is nonnegative, we obtain a
qualitatively sharp Gaussian upper bound. This is presented in
Section 5.

We will use the following notations throughout the paper. ${\bf
M}$ denotes a compact Riemannian manifold without boundary, unless
stated otherwise; $g, R_{ij}$ will be the metric and Ricci
curvature; $\nabla$, $\Delta$ the corresponding gradient and
Laplace-Beltrami operator; $c$ with or without index denote
generic positive constant that may change from line to line. In
case the metric $g(t)$ evolves with time, then $d(x, y, t)$ will
denote the corresponding distance function; $d\mu(x, t)$ denotes
the volume element under $g(t)$; We will still use $\nabla$,
$\Delta$ the corresponding gradient and Laplace-Beltrami operator,
when no confusion arises.

To close the introduction, we point out that all results in the
paper are stated for compact manifolds or bounded domains. However
similar results can be proven for noncompact manifolds under
appropriate assumptions near infinity.

\medskip

\section{log derivative estimates for Dirichlet
  heat kernel and Poisson heat kernel on domains.}

There have been several log gradient estimates available for the
heat kernel on complete manifolds, compact manifolds without
boundary and for the Neumann heat kernel. In the introduction, we
mentioned the papers [LY] and [Wa]. For compact manifolds without
boundary, we refer the reader to the papers [Sh], [H] (Corollary
1.3), [No], [MS], [Hs] and [ST]. However, an estimate for the
Dirichlet heat kernel is clearly missing. This is done in the next
theorem. The estimate is sharp in general as can be seen from the
heat kernel formula for the Euclidean half space. Let us mention
that for bounded domains, the large time behavior of heat kernels
is determined by the first eigenvalue and eigenfunction. So we
will only deal with the most interesting, small time case.

\begin{theorem}
  Let $D$ be a bounded $C^2$ domain in a Riemainnian manifold and
  $G=G(x, t; y, 0)$ and $P=P(x, t; y, 0)$ be the Dirichlet heat
  kernel and Poisson heat kernel respectively. Also let $\rho(x) = dist(x,
  \partial D)$ and $d(x, y)$ be the Riemannian distance. Given $T>0$, there exists
  a constant $C$ depending on $T$ and $D$ such that
  \[
  |\nabla_x \log G(x, t; y, 0)| \le
  \begin{cases} \frac{C}{\rho(x)}, \qquad \rho(x) \le \sqrt{t};\\
  \frac{C}{\sqrt{t}} \ [ 1 + \frac{d(x, y)}{\sqrt{t}} ], \qquad
  \rho(x) > \sqrt{t};
  \end{cases}
  \leqno(2.1)
  \]for all $x, y \in D$ and $0<t<T$; and for all $x \in D$, $y \in
  \partial D$ and $0<t<T$,
 \[
  |\nabla_x \log P(x, t; y, 0)| \le
  \begin{cases} \frac{C}{\rho(x)}, \qquad \rho(x) \le \sqrt{t};\\
  \frac{C}{\sqrt{t}} \ [ 1 + \frac{d(x, y)}{\sqrt{t}} ], \qquad
  \rho(x) > \sqrt{t}.
  \end{cases}
   \leqno(2.2)
  \]
  \proof
  \end{theorem}

Let us prove (2.1) first. The proof of (2.2) is similar and will
be sketched later. As explained earlier,the most interesting case
for the derivative estimate is for small time. Hence we can take
$T$ to be sufficiently small. Here we will take $T$ so small that
the boundary Harnack principle of [FGS] holds when $\rho(x) \le 2
T$. Here we notice that even though the boundary Harnack principle
was proven in the Euclidean case in that paper, it is still valid
in the current case. This is so because we can cover the boundary
of $D$ by a finite number of metric balls with radius less than
the injectivity radius. And then we can convert the
Laplace-Beltrami operator into an elliptic equation with smooth
coefficients in ${\bf R}^n$.

For a fixed $t_0 \in (0, T)$ and $y \in D$, we write
\[
f(x, t) = G(x, t; y, -t_0), \qquad x \in D, t>0;
\]
\[
\Omega_{t_0} = \{ (z, \tau) \ | \ x \in D, 0< \tau \le t_0, \
\rho(z) \ge \sqrt{\tau} \ \}. \leqno(2.3)
\]Fixing $(x, t) \in D \times [0, t_0] - \Omega_{t_0}$, we can apply
the gradient estimate in Theorem 1.1 of [SZ] on the cube
\[
Q_{x, t} = B(x, \rho(x)) \times [t-\frac{\rho^2(x)}{2}, t] \subset D
\times [-t_0, t_0].
\] This gives us
\[
  \frac{|\nabla f(x, t) |}{f(x, t)} \le \frac{C}{\rho(x)}
  \big{(} 1+ \log \frac{A}{f(x, t)} \big{)}
  \leqno(2.4)
\]Here $A = sup_{Q_{x, t}} f$. For a proof of (2.4)
 and that of  Theorem 1.1 in [SZ], please go to
Theorem 3.1  in the next section, which contains Theorem 1.1 in
[SZ] as a special case.

Now we apply the standard Harnack inequality of [LY] on manifold
to reach
\[
A = sup_{Q_{x, t}} f \le c_1 f(x, t+ \rho(x)^2).
\]Then the boundary Harnack inequality of [FGS] gives us
\[
 f(x, t+ \rho(x)^2) \le c_2 f(x, t)
 \]since $f$ vanishes on $\partial D \times (-t_0, t_0)$. Therefore
 \[
 A = sup_{Q_{x, t}} f \le c_3 f(x, t)
\leqno(2.5)
\]Substituting (2.5) to (2.4), we deduce,
for $(x, t) \in D \times [0, t_0] - \Omega_{t_0}$,
\[
  \frac{|\nabla f(x, t) |}{f(x, t)} \le \frac{C}{\rho(x)}
  \leqno(2.6)
\]This proves the first part of (2.1).

Next we work in $\Omega_{t_0}$.  Let us observe that on the sides of
$\partial \Omega_{t_0}$, there holds $\rho(x) = \sqrt{t}$. Hence for
such $x$, inequality (2.6) becomes
\[
  \frac{|\nabla f(x, t) |}{f(x, t)} \le \frac{C}{\sqrt{t}}
\]i.e.
\[
 \frac{|\nabla f(x, t) |^2}{f(x, t)} \le \frac{C}{t} f(x, t), \qquad
 \rho(x) = \sqrt{t}.
 \leqno(2.7)
\]Let $m = \sup_{\Omega_{t_0}} f$, then for any $b>0$, we have
\[
f \log \frac{b m}{f} \ge f \log b.
\]Now we use the calculation in the proof of Theorem 1.1 in [H]
(p115) to reach
\[
\begin{aligned}
 &\Delta (f \log \frac{b m}{f}) - \partial_t (f \log \frac{b
m}{f}) \\
&=(\Delta f -\partial_t f) \log b + (\Delta -\partial_t) (f
\log\frac{m}{f})\\
&=-\frac{|\nabla f|^2}{f}.
 \end{aligned}
 \leqno(2.9)
\]Also
\[
\al (\Delta -\partial_t)&(\frac{|\nabla f|^2}{f}) = \frac{2}{f}
\bigg{|}\partial_i \partial_j f - \frac{\partial_i f \partial_j
f}{f} \bigg{|}^2 + 2 R_{ij} \frac{\partial_i f \partial_j f}{f}\\
& \ge - 2 K \frac{|\nabla f|^2}{f}.
 \eal
 \leqno(2.10)
\]Here $-K$ is the lower bound of the Ricci curvature.
Therefore, for
\[
h = \frac{t}{1+ 2K t} \frac{|\nabla f|^2}{f} - f \log \frac{b m}{f},
\]we have
\[
\Delta h(x, t) - \partial_t h(x, t) \ge 0, \qquad (x, t) \in
\Omega_{x_0}.
\]We $t =0$, it is clear that $h \le 0$. On the sides of $\partial
\Omega_{t_0}$, i.e., when $\rho(x) = \sqrt{t}$, one can choose $b$
sufficiently large so that
\[
f \log \frac{b m}{f} \ge C f \ge t \frac{|\nabla f|^2}{f}.
\]Here we just used (2.7) and the constant $C$ is from there too.
Therefore $h \le 0$ on the sides of $\partial \Omega_{t_0}$. When
$T$ is sufficiently small, we know that $\partial \Omega_{t_0}$ is
connected and we can apply the maximum principle on this time
dependent domain to conclude that $h(x) \le 0$ in $\Omega_{t_0}$.
i.e.
\[
t \frac{|\nabla f|^2}{f} \le C(T, K) f \log \frac{b m}{f}.
\leqno(2.11)
\]In particular, this shows, with
\[
t=t_0,  f=G(x, t_0; y, -t_0) = G(x, 2 t_0; y, 0),
\]that
\[
\frac{|\nabla_x G(x, 2t_0; y, 0)|}{G(x, 2t_0; y, 0)} \le \frac{C(T,
K)}{\sqrt{t_0}} \bigg{[} 1 + \sqrt{\log \frac{m}{G(x, 2t_0; y, 0)}}
\bigg{]}
\]when $\rho(x) \ge \sqrt{2 t_0}$. Here
\[
m = \sup_{\Omega_{t_0}} f = \sup_{\rho(z) \ge \sqrt{\tau}, 0<\tau
\le t_0} f(z, \tau) =  \sup_{\rho(z) \ge \sqrt{\tau}, 0<\tau \le
t_0} G(z, \tau; y, -t_0).
\]Making a change of variables $t_0 \to t_0/2$, we have
\[
\frac{|\nabla_x G(x, t_0; y, 0)|}{G(x, t_0; y, 0)} \le \frac{C(T,
K)}{\sqrt{t_0}} \bigg{[} 1 + \sqrt{\log \frac{m}{G(x, t_0; y, 0)}}
\bigg{]} \leqno(2.12)
\]when $\rho(x) \ge \sqrt{t_0}$.

By the Dirichlet heat kernel upper bound in Davies [Da], we know
that
\[
m \le C \bigg{(} \frac{\rho(y)}{\sqrt{t_0}} \wedge 1  \bigg{)}
\frac{1}{|B(y, \sqrt{t_0})|}. \leqno(2.13)
\]By the lower bound estimate in [Z], there holds
\[
G(x, t_0; y, 0) \ge C \bigg{(} \frac{\rho(y)}{\sqrt{t_0}} \wedge 1
\bigg{)} \frac{1}{|B(y, \sqrt{t_0})|} e^{- c d(x, y)^2/t_0}.
\leqno(2.14)
\]Here we note that the lower bound was proven under the assumption
that the Ricci curvature is nonnegative. However for short time
behavior this assumption is not necessary. Substituting (2.13) and
(2.14) to (2.12), we obtain
\[
\frac{|\nabla_x G(x, t_0; y, 0)|}{G(x, t_0; y, 0)} \le \frac{C(T,
K)}{\sqrt{t_0}} \bigg{[} 1 + \frac{d(x, y)}{\sqrt{t_0}} \bigg{]}
\leqno(2.15)
\]when $\rho(x) \ge \sqrt{t_0}$. Now (2.1) follows from (2.7) and
(2.15).

To prove (2.2), let us recall the results in [Da] (upper bound) and
[Z] (lower bound): there exists $c_1$ and $c_2$ such that
\[
\al
& \frac{1}{c_1} \bigg{(} \frac{\rho(x)}{\sqrt{t}} \wedge 1
\bigg{)} \bigg{(} \frac{\rho(y)}{\sqrt{t}} \wedge 1 \bigg{)}
\frac{1}{|B(y, \sqrt{t})|} e^{-  d(x, y)^2/(c_2 t)} \le \\
&\qquad G(x, t; y, 0) \le c_1 \bigg{(} \frac{\rho(x)}{\sqrt{t}}
\wedge 1 \bigg{)} \bigg{(} \frac{\rho(y)}{\sqrt{t}} \wedge 1
\bigg{)} \frac{1}{|B(y, \sqrt{t})|} e^{- c_2 d(x, y)^2/t}. \eal
\]for all $x, y \in D$ and $0<t \le T$.

Given $y \in \partial D$, the Poisson heat kernel is defined as
\[
P(x, t; y, s) = -\frac{\partial}{\partial_{n_y}} G(x, t; y, 0).
\]Therefore one has the two-sided bound
\[
\frac{1}{c_1} \bigg{(} \frac{\rho(x)}{\sqrt{t}} \wedge 1 \bigg{)}
\frac{1}{|B(y, \sqrt{t})|} e^{-  d(x, y)^2/(c_2 t)} \le P(x, t; y,
0) \le c_1 \bigg{(} \frac{\rho(x)}{\sqrt{t}} \wedge 1 \bigg{)}
 \frac{1}{|B(y,
\sqrt{t})|} e^{- c_2 d(x, y)^2/t}.
\]The rest of the proof for (2.2) is identical to that of (2.1).
\qed \bigskip

Our next theorem provides a sharpened local Li-Yau estimate. In
1986 Li and Yau proved the following famous estimate.
\medskip

 {\it {\bf Theorem }(Li-Yau [LY]). Let ${\bf M}$ be a complete
manifold with dimension $n \ge 2$, $Ricci ({\bf M}) \ge -K$, $K
\ge 0$. Suppose $u$ is any positive solution to the heat equation
in $B(x_0, R) \times [t_0-T, t_0] \subset {\bf M} \times [t_0-T,
t_0] $. Then, for any $\alpha \in (0, 1)$, there exists a constant
$c=c(n, \alpha)$ such that
\[
\alpha \frac{|\nabla u|^2}{u^2} -\frac{u_t}{u} \le \frac{c}{R^2}+
\frac{c}{T}+ c K,
\]in $B(x_0, R/2) \times [t_0-T/2, t_0]$.

Moreover, if ${\bf M}$ has nonnegative Ricci curvature and
$R=\infty$, i.e. $B(x_0, R) = {\bf M}$, then
\[
 \frac{|\nabla u|^2}{u^2} -\frac{u_t}{u} \le
\frac{c_n}{T}.
\]}

\medskip

Let us observe that, even in the case of nonnegative Ricci
curvature, the first local estimate does not match the second global
estimate completely, due to the presence of the parameter
$\alpha<1$. Here we show that $\alpha$ can be taken as $1$ modulo a
lower order term.  We mention that our estimate in the next theorem
is new only in the local sense.  The global estimate was already
proven in [Y] by a using a more involved quantity. The very short
proof, simpler than previous ones, is based on a modification of an
idea in [H] and the cut-off method in [LY].

\begin{theorem}
Let $B(x_0, R)$ be a geodesic ball in a Riemannian  manifold ${\bf
M}$ with dimension $n \ge 2$ such that $Ricci |_{B(x_0, R)} \ge
-K$, $K \ge 0$. Suppose $u$ is any positive solution to the heat
equation in $B(x_0, R) \times [t_0-T, t_0]$. Then
\[
\frac{|\nabla u|^2}{u^2} - \frac{u_t}{u} \le \frac{c_n}{R^2} +
\frac{c_n}{T} + c_n K + c_n \sqrt{K} \sup \frac{|\nabla u|}{u} +
\frac{c_n}{R} \sup \frac{|\nabla u|}{u},
\]
in $B(x_0, R/2) \times [t_0-T/2, t_0]$. Here $c_n$ depends only on
the dimension $n$. \proof
\end{theorem}

 By direct computation (see [H]), we have
\[
(\Delta -\partial_t)(\frac{|\nabla u|^2}{u}) = \frac{2}{u}
\bigg{|}\partial_i \partial_j u - \frac{\partial_i u \partial_j
u}{u} \bigg{|}^2 + 2 R_{ij} \frac{\partial_i u \partial_j u}{u}.
\]In view of the estimate
\[
\bigg{|}\partial_i \partial_j u - \frac{\partial_i u \partial_j
u}{u} \bigg{|}^2 \ge \frac{1}{n} \bigg{(} \Delta u - \frac{|\nabla
u|^2}{u} \bigg{)}^2,
\]the above implies
\[
 (\Delta -\partial_t)(\frac{|\nabla u|^2}{u}) \ge \frac{2}{n u}
 \bigg{(} \Delta u - \frac{|\nabla u|^2}{u} \bigg{)}^2 +
 2 R_{ij} \frac{\partial_i u \partial_j
u}{u}.
\]Since $\Delta u$ is also a solution to the heat equation,
it follows that
\[
(\Delta -\partial_t)(- \Delta u + \frac{|\nabla u|^2}{u}) \ge
\frac{2}{n u}
 \bigg{(} \Delta u - \frac{|\nabla u|^2}{u} \bigg{)}^2  -
  2 K \frac{|\nabla u|^2}{u}.
\]Let us write
\[
q=- \Delta u + \frac{|\nabla u|^2}{u} = \frac{|\nabla u|^2}{u} -
u_t.
\]Then $q$ satisfies
\[
(\Delta -\partial_t) q \ge  \frac{2}{n u} q^2  - 2 K \frac{|\nabla
u|^2}{u}.
\]Define
\[
H= q/u.
\]Then $H$ satisfies
\[
(\Delta -\partial_t) H \ge  \frac{2}{n} H^2  - 2 K \frac{|\nabla
u|^2}{u^2} - 2 \nabla H \nabla \ln u. \leqno(2.16)
\]

Now we can use the Li-Yau idea of cut-off functions to derive the
desired bound. The only place that may cause difficulty is that $H$
may change sign. However it turns out that it does not hurt. Here is
the detail. Let $\psi=\psi(x, t)$ be a smooth cut-off function
supported in $Q_{R, T} \equiv B(x_0, R) \times [t_0-T, t_0]$,
satisfying the following properties

(1). $\psi = \psi(d(x, x_0), t) \equiv \psi(r, t)$; $\psi(x, t) =
1$ in $Q_{R/2, T/4}$, \ $0 \le \psi \le 1$.

(2). $\psi$ is decreasing as a radial function in the spatial
variables.

(3). $\frac{|\partial_r \psi|}{\psi^{a}} \le \frac{C_a}{R}$,
$\frac{|
\partial^2_r \psi|}{\psi^{a}} \le \frac{C_a}{R^2}$when $0<a<1$.

(4). $\frac{|\partial_t \psi|}{\psi^{1/2}} \le \frac{C}{T}$.

  Then, from
(2.16) and a straight forward calculation, one has
\[
\aligned
  \Delta &(\psi H)  - (\psi H)_t - 2 \frac{\nabla
  \psi}{\psi}
  \cdot  \nabla (\psi H) + 2 \psi K \frac{|\nabla u|^2}{u^2}+ 2 \nabla(\psi H)
  \nabla \ln u\\
&\ge \frac{2}{n} \psi H^2 + (\Delta \psi) H - 2 \frac{|\nabla
\psi|^2}{\psi}
H - \psi_t H + 2 H \nabla \psi \nabla \ln u\\
&=\frac{2}{n} \psi H^2 - 2 \frac{|\nabla \psi|^2}{\psi} H +
(\partial^2_r \psi + (n-1) \frac{
\partial_r \psi}{r} + \partial_r \psi \partial_r \log \sqrt{g} ) H
- \psi_t H + 2 H \nabla \psi \nabla \ln u.
\endaligned
\leqno(2.17)
\]Suppose that at $(y, s)$, the function $\psi H$  reaches a maximum.
If the value is non-positive, there is nothing to prove. So we
assume the maximum value is positive. Then (2.17) shows
\[
2 \psi K \frac{|\nabla u|^2}{u^2} + 2 \frac{|\nabla \psi|^2}{\psi} H
\ge \frac{2}{n} \psi H^2 + (\partial^2_r \psi + (n-1) \frac{
\partial_r \psi}{r} + \partial_r \psi \partial_r \log \sqrt{g} ) H
- \psi_t H +  2 H \nabla \psi \nabla \ln u.
\]
In the above, the only term we need extra care of is
\[
\partial_r \psi \partial_r \log \sqrt{g}  H.
\]Note that $-C/R \le \partial_r \psi /\psi^a \le 0$, $\partial_r \log \sqrt{g} \le
\sqrt{K}$ and $H(y, s)>0$. Therefore
\[
\al
 2 \psi K \frac{|\nabla u|^2}{u^2} &+  2 \frac{|\nabla
\psi|^2}{\psi} H + 2 \sqrt{\psi} H \frac{\nabla \psi}{\sqrt{\psi}}
\nabla \ln u \\
&\ge \frac{2}{n u} \psi H^2 + (\partial^2_r \psi + (n-1) \frac{
\partial_r \psi}{r}) H - C \sqrt{K} \psi^a H/R
- \psi_t H. \eal
\]This shows that
\[
\psi H^2 = \psi \bigg{(} \frac{|\nabla u|^2}{u^2} - \frac{u_t}{u}
\bigg{)}^2 \le (\frac{c_n}{R^4} + \frac{c_n}{T^2} + c_n K^2)  + c_n
K \frac{|\nabla u|^2}{u} + \big{(} \frac{c_n}{R} |\nabla \ln u|
\big{)}^2.
\]Hence
\[
\frac{|\nabla u|^2}{u^2} - \frac{u_t}{u} \le \frac{c_n}{R^2} +
\frac{c_n}{T} + c_n K + c_n \sqrt{K} \frac{|\nabla u|}{u} +
\frac{c_n}{R} \frac{|\nabla u|}{u}.
\] in the half parabolic cube.
\qed

\section{gradient estimates on the log
 temperature under backward and
forward Ricci flow}

In this section we will prove certain localized or global gradient
bound on the heat equation under backward and forward Ricci flow,
i.e. equations (3.1) and (3.2) below. This estimate is a
generalization of the results in [H] and [SZ], where the heat
equation under a fixed metric is studied. Similar estimates for the
conjugate heat equation (i.e. when $\Delta$ is replaced by $\Delta
-R$ in (3.1) or (3.2)) were proven in [Ni3], [CKNT] and [CCGGIIKLLN]
Chapter 8. This estimate then also relies on the derivative of the
scalar curvature $R$.

The current estimate under the forward Ricci flow ((3.2)) will be
useful for Section 5, where we will prove a global Gaussian upper
estimate for Perelman's equation under nonnegative Ricci curvature
assumption.

Recall that the heat equation under backward and forward Ricci
flow are given by
\[
\begin{cases}
\Delta u - \partial_t u = 0, \\
\frac{d}{dt} g_{ij} =  2 R_{ij}
\end{cases}
\leqno(3.1)
\]and
\[
\begin{cases}
\Delta u - \partial_t u = 0, \\
\frac{d}{dt} g_{ij} =  - 2 R_{ij}.
\end{cases}
\leqno(3.2)
\]

For (3.1) we have the following:

\medskip

\begin{theorem}
Let ${\bf M}$ be a compact Riemannian manifold equipped with a
family of Riemannian metric evolving under the backward Ricci flow
in (3.1).

 (a) (local estimate) . Suppose $u$ is any positive
solution to (3.1)  in
\[
Q_{R, T}= \{ (x, t) \ | \ x \in {\bf M}, d(x, x_0, t)< R,  t \in
[t_0-T, t_0]\}
\]such that the $Ricci \ge -k$ throughout.
   Suppose also $u \le M$ in  $Q_{R, T}$. Then there exists a dimensional
  constant $c$ such that
\[
  \frac{|\nabla u(x, t) |}{u(x, t)} \le c (\frac{1}{R} +
\frac{1}{T^{1/2}} +\sqrt{k} )
   \big{(} 1+ \log \frac{M}{u(x, t)} \big{)}
  \]in $Q_{R/2, T/2}$.

  (b). (global estimate) Suppose $u$ is any positive
solution to (3.1)  in ${\bf M} \times [0, T]$. Under the assumption
that $Ricci \ge 0$, it holds
\[
\frac{|\nabla u(x, t) |}{u(x, t)} \le   \frac{1}{t^{1/2}}
    \sqrt{\log \frac{M}{u(x, t)}}
\]for $M = \sup_{{\bf M} \times [0, T]} u$ and $(x, t) \in {\bf M} \times [0,
T]$.
\end{theorem}

{\it Remark.} As pointed out in [SZ], the local and global estimate
can not replace each other. Also note that there is no other
curvature assumption in part (b), nor any constants.
\medskip

{\bf Proof of Theorem 3.1 (a).}

We will use the idea in [SZ] with certain modifications to handle
the changing nature of the metric.
  Suppose
$u$ is a solution to the heat equation in the statement of the
theorem in the parabolic cube $ Q_{R, T}$.
It is clear that the
gradient estimate in Theorem 3.1 is invariant under the scaling $u
\to u/M$. Therefore, we can and do assume that $0< u \le 1$.

Write
\[
f = \log u, \qquad w \equiv | \nabla \log (1-f)|^2 = \frac{|\nabla
f|^2}{(1-f)^2}.
\]

Since $u$ is a solution to the heat equation, simple calculation
shows that
\[
\Delta f + | \nabla f|^2 - f_t =0.
\]

We will derive an equation for $w$. First notice that
\[
\al w_t &= \frac{2 \nabla f (\nabla f)_t}{(1-f)^2} + \frac{2
|\nabla f|^2 f_t}{(1-f)^3} + \frac{2 Ric (\nabla f, \nabla f)}{(1-f)^2}\\
&= \frac{2 \nabla f \nabla (\Delta f + |\nabla f|^2)}{(1-f)^2} +
\frac{2 |\nabla f|^2 (\Delta f + |\nabla f|^2)}{(1-f)^3} - \frac{2
Ric (\nabla f, \nabla f)}{(1-f)^2} \eal
\]In local orthonormal system, this can be written as
\[
w_t=\frac{2  f_j f_{iij} + 4 f_i f_j f_{ij}}{(1-f)^2} + 2 \frac{
f^2_i f_{jj} + |\nabla f|^4}{(1-f)^3} - \frac{2 R_{ij} f_i f_j
}{(1-f)^2}. \leqno(3.3)
\]Here and below, we have adopted the convention $f^2_i = |\nabla
f|^2$ and $f_{ii} = \Delta f$.

Next
\[
\nabla w = \big{(} \frac{f^2_i}{(1-f)^2} \big{)}_j = \frac{2 f_i
f_{ij}}{(1-f)^2} + 2 \frac{f^2_i f_j}{(1-f)^3}. \leqno(3.4)
\]It follows that
\[
\al \Delta w &= \big{(} \frac{f^2_i}{(1-f)^2} \big{)}_{jj}\\
&= \frac{ 2 f^2_{ij}}{(1-f)^2} + \frac{2 f_i f_{ijj}}{(1-f)^2} +
\frac{4 f_i f_{ij} f_j}{(1-f)^3}\\
&\qquad + \frac{4 f_i f_{ij} f_j}{(1-f)^3} + 2 \frac{f^2_i
f_{jj}}{(1-f)^3} + 6 \frac{f^2_i f^2_j}{(1-f)^4}. \eal \leqno(3.5)
\]By (3.5) and (3.3),
\[
\al
&\Delta w - w_t \\
&=\frac{ 2 f^2_{ij}}{(1-f)^2} + 2 \frac{ f_i f_{ijj}-f_j
f_{iij}}{(1-f)^2}\\
&\qquad + 6 \frac{|\nabla f|^4}{(1-f)^4} + 8 \frac{ f_i f_{ij}
f_j}{(1-f)^3} + 2 \frac{f^2_i f_{jj}}{(1-f)^3} \\
&\qquad - 4 \frac{ f_i f_{ij} f_j}{(1-f)^2} - 2 \frac{f^2_i
f_{jj}}{(1-f)^3} - 2 \frac{|\nabla f|^4}{(1-f)^3} + \frac{2 R_{ij}
f_i f_j }{(1-f)^2}. \eal
\]The 5th and 7th terms on the righthand side of this identity
cancel each other. Also, by Bochner's identity
\[
f_i f_{ijj}-f_j f_{iij} = f_j( f_{jii}-f_{iij}) = R_{ij} f_i f_j.
\]So the second term  doubles with the last term.
Therefore
\[
\al
&\Delta w - w_t \\
& = \frac{ 2 f^2_{ij}}{(1-f)^2} + 6 \frac{|\nabla f|^4}{(1-f)^4} + 8
\frac{ f_i f_{ij} f_j}{(1-f)^3} - 4 \frac{ f_i f_{ij} f_j}{(1-f)^2}
- 2 \frac{|\nabla f|^4}{(1-f)^3} + \frac{4 R_{ij} f_i f_j
}{(1-f)^2}. \eal \leqno(3.6)
\]

Notice from (3.4) that
\[
\nabla f \nabla w = \frac{2 f_i f_{ij}f_j }{(1-f)^2} + 2
\frac{f^2_i f^2_j}{(1-f)^3}.
\]Hence
\[
0 = 4 \frac{ f_i f_{ij} f_j}{(1-f)^2}  - 2 \nabla f \nabla w + 4
\frac{|\nabla f|^4}{(1-f)^3}, \leqno(3.7)
\]
\[
0= -4 \frac{ f_i f_{ij} f_j}{(1-f)^3} + [ 2 \nabla f \nabla w - 4
\frac{|\nabla f|^4}{(1-f)^3} ] \frac{1}{1-f}. \leqno(3.8)
\]Adding (3.6) with (3.7) and (3.8), we deduce
\[
\al
&\Delta w - w_t \\
& = \frac{ 2 f^2_{ij}}{(1-f)^2} + 2 \frac{|\nabla f|^4}{(1-f)^4}
+ 4 \frac{ f_i f_{ij} f_j}{(1-f)^3}  \\
&\qquad + \frac{2}{1-f} \nabla f \nabla w - 2 \nabla f \nabla w + 2
\frac{|\nabla f|^4}{(1-f)^3} + \frac{4 R_{ij} f_i f_j }{(1-f)^2}.
\eal
\]Since
\[
\frac{ 2 f^2_{ij}}{(1-f)^2} + 2 \frac{|\nabla f|^4}{(1-f)^4} + 4
\frac{ f_i f_{ij} f_j}{(1-f)^3} \ge 0,
\]we have
\[
\Delta w - w_t
  \ge  \frac{2 f}{1-f} \nabla f \nabla w  +
2 \frac{|\nabla f|^4}{(1-f)^3} - 4 k w.
\]Since $f \le 0$, it follows that
\[
\Delta w - w_t
  \ge  \frac{2 f}{1-f} \nabla f \nabla w  +
2 (1-f) \frac{|\nabla f|^4}{(1-f)^4} - 4 k w,
\]i.e.
\[
\Delta w - w_t
  \ge  \frac{2 f}{1-f} \nabla f \nabla w  +
2 (1-f) w^2 - 4 k w. \leqno(3.9)
\]

From here, we will use a cut-off function  to derive the desired
bounds.
  Let $\psi=\psi(x, t)$ be a smooth cut-off function
supported in $Q_{R, T}$, satisfying the following properties

(1). $\psi = \psi(d(x, x_0, t), t) \equiv \psi(r, t)$; $\psi(x, t)
= 1$ in $Q_{R/2, T/4}$, \ $0 \le \psi \le 1$.

(2). $\psi$ is decreasing as a radial function in the spatial
variables.

(3). $\frac{|\partial_r \psi|}{\psi^{a}} \le \frac{C_a}{R}$,
$\frac{|
\partial^2_r \psi|}{\psi^{a}} \le \frac{C_a}{R^2}$when $0<a<1$.

(4). $\frac{|\partial_t \psi|}{\psi^{1/2}} \le \frac{C}{T}$.

  Then, from
(3.9) and a straight forward calculation, one has
\[
\aligned
  \Delta (\psi w)& + b \cdot \nabla (\psi w) - 2 \frac{\nabla
  \psi}{\psi}
  \cdot  \nabla (\psi w) - (\psi w)_t\\
&\ge 2 \psi (1-f) w^2 + (b  \cdot \nabla \psi) w - 2 \frac{|\nabla
\psi|^2}{\psi} w + (\Delta \psi) w - \psi_t w -  4 k \psi w,
\endaligned
\]where we have written
\[
b = - \frac{2 f}{1-f} \nabla f.
\]Comparing with the heat equation under a fixed metric,
 the last term $-\psi_t w$ is more complicated. It is given by
 \[
-\psi_t w = -[ \frac{\partial \psi }{\partial t}  + \frac{\partial
\psi}{\partial r}  \  \frac{\partial d(x, x_0, t)}{\partial t} ] \
w.
\]By our assumption that $ Ricci \ge -k$ and that $-c/R \le
\frac{\partial \psi}{\partial r} \le 0$, we
have
\[
-\psi_t w \ge - \frac{\partial \psi }{\partial t} w - c k w
\psi^{1/2}.
\]Here $\frac{\partial \psi }{\partial t} \equiv
 \frac{\partial \psi(r, t) }{\partial t}$.

Therefore
\[
\aligned
  \Delta (\psi w)& + b \cdot \nabla (\psi w) - 2 \frac{\nabla
  \psi}{\psi}
  \cdot  \nabla (\psi w) - (\psi w)_t\\
&\ge 2 \psi (1-f) w^2 + (b  \cdot \nabla \psi) w - 2 \frac{|\nabla
\psi|^2}{\psi} w + (\Delta \psi) w - \frac{\partial \psi }{\partial
t} w - c k w \psi^{1/2}.
\endaligned
\leqno(3.10)
\]

Suppose the  maximum of $\psi w$ is reached at $(x_1, t_1)$. By
[LY], we can assume, without loss of generality that $x_1$ is not
in the cut-locus of ${\bf M}$. Then at this point, one has,
  $\Delta (\psi w) \le 0$, $(\psi w)_t \ge 0$ and $\nabla (\psi
  w)=0$. Therefore
  \[
2 \psi (1-f) w^2(x_1, t_1) \le - [ \ (b  \cdot \nabla \psi) w - 2
\frac{|\nabla \psi|^2}{\psi} w + (\Delta \psi) w - \frac{\partial
\psi }{\partial t} w \ ](x_1, t_1) + c k w \psi^{1/2}. \leqno(3.11)
\]We need to find an upper bound for each term of the righthand
side of (3.11).

\[
\al |(b  \cdot \nabla \psi) w| &\le \frac{ 2|f|}{1-f} |\nabla f| w
|\nabla \psi| \le 2 w^{3/2} |f| \  |\nabla \psi| \\
&= 2 [\psi (1-f) w^2]^{3/4} \ \frac{f |\nabla \psi|}{[\psi
(1-f)]^{3/4}}\\
&\le \psi (1-f) w^2 + c \frac{(f |\nabla \psi|)^4}{[\psi
(1-f)]^3}. \eal
\]This implies
\[
|(b  \cdot \nabla \psi) w| \le (1-f) \psi w^2 + c \frac{f^4}{R^4
(1-f)^3}. \leqno(3.12)
\]

For the second term on the righthand side of (3.11), we proceed as
follows
\[
\aligned
  \frac{|\nabla \psi|^2}{\psi} w &= \psi^{1/2} w
\frac{|\nabla
\psi|^2}{\psi^{3/2}}\\
&\le \frac{1}{8} \psi w^2 + c \big{(} \frac{|\nabla
\psi|^2}{\psi^{3/2}} \big{)}^2 \le \frac{1}{8} \psi w^2  + c
\frac{1}{R^4}.
\endaligned
\leqno(3.13)
\]

Furthermore, by the properties of $\psi$ and the assumption of on
the Ricci curvature, one has
\[
\aligned - (\Delta \psi) w &= -(\partial^2_r \psi + (n-1) \frac{
\partial_r \psi}{r} + \partial_r \psi \partial_r \log \sqrt{g} ) w\\
&\le (|\partial^2_r \psi| + 2 (n-1) \frac{|
\partial_r \psi|}{R}  ) w + \frac{c}{R} \sqrt{k} w \sqrt{\psi} \\
  &\le \psi^{1/2} w \frac{|\partial^2_r
\psi|}{\psi^{1/2}} + \psi^{1/2} w 2 (n-1) \frac{|
\partial_r \psi|}{R \psi^{1/2}} + \frac{c}{R} \sqrt{k} w \sqrt{\psi} \\
&\le \frac{1}{8} \psi w^2  + c \big{(} [\frac{|\partial^2_r
\psi|}{\psi^{1/2}}]^2 + [ \frac{|
\partial_r \psi|}{R \psi^{1/2}}]^2 + \frac{ c k}{R^2}.
  \endaligned
\]Therefore
\[
- (\Delta \psi) w \le \frac{1}{8} \psi w^2  + c \frac{1}{R^4}.
\leqno(3.14)
\]

Now we estimate $|\frac{\partial \psi }{\partial t}| \  w$.
\[
\aligned
  |\frac{\partial \psi }{\partial t}| \ w &= \psi^{1/2} w \frac{|\frac{\partial \psi }{\partial t}|}{\psi^{1/2}}\\
&\le \frac{1}{8} \big{(} \psi^{1/2} w \big{)}^{2} + c \big{(}
\frac{|\frac{\partial \psi }{\partial t}|}{\psi^{1/2}}
\big{)}^{2}.
\endaligned
\]This shows
\[
|\frac{\partial \psi }{\partial t}| w  \le \frac{1}{8} \psi w^2 +
c \frac{1}{T^2}. \leqno(3.15)
\]

Substituting (3.12)-(3.15) to the righthand side of (3.11), we
deduce,
\[
  2 (1-f) \psi w^2 \le (1-f)  \psi w^2 + c \frac{f^4}{R^4
(1-f)^3} + \frac{1}{2} \psi w^2 + \frac{c}{R^4} + \frac{c}{T^2} +
\frac{c k}{R^2} + k w \sqrt{\psi}.
\]Recall that $f \le 0$, therefore the above implies
\[
  \psi w^2(x_1, t_1) \le  c \frac{f^4}{R^4 (1-f)^4} +
   \frac{1}{2} \psi w^2(x_1, t_1) +
\frac{c}{R^4} +  \frac{c}{T^2} + c k^2.
\]Since $\frac{f^4}{(1-f)^4} \le 1$, the above shows, for all $(x,
t)$ in $Q_{R, T}$,
\[
\al \psi^2(x, t) w^2(x, t) &\le \psi^2(x_1, t_1) w^2(x_1, t_1)
\\
&\le \psi(x_1, t_1) w^2(x_1, t_1) \\
&\le c \frac{c}{R^4} + \frac{c}{T^2} + c k^2. \eal
\]Notice that $\psi(x, t) = 1$ in $Q_{R/2, T/4}$ and $w =|\nabla
f|^2/(1-f)^2$. We finally have
\[
\frac{|\nabla f(x, t)|}{1-f(x, t)} \le
  \frac{c}{R} + \frac{c}{\sqrt{T}} + c \sqrt{k}.
\]We have completed the proof of Theorem 3.1 (a) since $f=\log (u/M)$ with $M$
scaled to $1$.
\medskip

{\bf Proof of Theorem 3.1 (b).}

The proof is almost identical to that of Theorem 1.1 in [H] except
for an additional curvature term. By direct computation, we have
\[
(\Delta -\partial_t)(\frac{|\nabla u|^2}{u}) \ge \frac{2}{u}
\bigg{|}\partial_i \partial_j u - \frac{\partial_i u \partial_j
u}{u} \bigg{|}^2.
\]In the above, comparing with the fixed curvature case (2.10),
there is no more term containing the Ricci curvature.  By (2.9),
it holds
\[
\Delta (u \log \frac{M}{u}) - \partial_t (u \log \frac{M}{u})
=-\frac{|\nabla u|^2}{u}.
\]Since
\[
(\Delta -\partial_t)(t \frac{|\nabla u|^2}{u}) \ge  -\frac{|\nabla
u|^2}{u}
\]the maximum principle implies that
\[
\frac{|\nabla u|^2}{u^2} \le \frac{1}{t} u \log \frac{M}{u}. \qed
\]

\bigskip

The remainder of the section deals with (3.2). For (3.2), we no
longer have the nice cancelation effect that associated with (3.1).
So we only obtain the following global gradient estimate under
curvature assumptions.

\medskip

\begin{theorem}
Let ${\bf M}$ be a complete Riemannian manifold equipped with a
family of Riemannian metric evolving under the forward Ricci flow in
(3.2) with $t \in [0, T]$. Suppose $u$ is any positive solution to
(3.2) in ${\bf M} \times [0, T]$. Then, it holds
\[
\frac{|\nabla u(x, t) |}{u(x, t)} \le  \sqrt{ \frac{1}{t}}
    \sqrt{\log \frac{M}{u(x, t)}}
\]for $M = \sup_{{\bf M} \times [0, T]} u$ and $(x, t) \in {\bf M} \times [0,
T]$.

Moreover,  the following interpolation inequality holds for any
$\delta>0$, $x, y \in {\bf M}$ and $0< t \le T$:
\[
u(y, t) \le c_1 u(x, t)^{1/(1+\delta)} M^{\delta/(1+\delta)} e^{c_2
d(x, y, t)^2/t}.
\]Here $c_1, c_2$ are positive constants depending only on
$\delta$.
\end{theorem}
\medskip

{\bf Proof of Theorem 3.2.}

This again is almost the same as that of Theorem 1.1 in [H]. By
direct calculation
\[
\Delta (u \log \frac{M}{u}) - \partial_t (u \log \frac{M}{u})
=-\frac{|\nabla u|^2}{u},
\]
\[
(\Delta -\partial_t)(\frac{|\nabla u|^2}{u}) = \frac{2}{u}
\bigg{|}\partial_i \partial_j u - \frac{\partial_i u \partial_j
u}{u} \bigg{|}^2 \\
 \ge 0.
\]The first inequality follows immediately from the maximum principle since
\[
t \frac{|\nabla u|^2}{u} - u \log \frac{M}{u}
\]is a sub-solution of the heat equation.

To prove the second inequality, we set
\[
l(x, t) = \log (M/u(x, t)).
\]Then the first inequality implies
\[
|\nabla \sqrt{l(x, t)} | \le 1/\sqrt{t}.
\]Fixing two points $x$ and $y$, we can integrate along a geodesic to
reach
\[
\sqrt{\log (M/u(x, t))} \le \sqrt{\log (M/u(y, t))} + \frac{d(x, y,
t)}{\sqrt{t}}.
\]The result follows by squaring both sides.
\qed

\section {Pointwise and gradient estimate for the fundamental solution
to an equation of  Perelman's}

In the paper [P] Perelman introduced an equation which after time
reversal becomes
\[
\begin{cases}
\Delta u - R u - \partial_t u = 0,\\
\frac{d}{dt} g_{ij} = 2 R_{ij}. \end{cases}
 \leqno(4.0)
\]Here as before $\Delta$ is the Laplace-Beltrami operator with
respect to the metric $g_{ij}$ evolving by the backward Ricci
flow. $R$ is the scalar curvature.  This equation and the
associated monotonicity formula have proven to be of fundamental
importance. Using the maximum principle and reduced distance,
Perelman proved a lower bound for the fundamental solution to
(4.0). An outstanding feature of the estimate is that it does need
any explicit curvature assumption. The information on curvature is
encoded in the reduced distance. From the analysis point of view,
it would be desirable to establish an upper bound for the
fundamental solution too. Here we first prove an upper bound under
no explicit curvature assumptions. The bound is in terms of  more
traditional geometric quantities, i.e. the best constant in
Sobolev imbedding or Yamabe constant, which are controlled by the
lower bound of the Ricci curvature and injectivity radius. Under
more restrictive curvature assumptions, we are able to prove a
Gaussian like upper bound. Let us mention the method by maximum
principle alone does not seem to yield the upper bound. Our method
is based on the one by J. Nash. For related results on local lower
and upper bounds for fundamental solutions of (3.2) and for a
global lower bound for the conjugate of (3.2) in the spirit of
Perelman, please see the interesting papers [G], [Ni1] and [Ni2].

In order to state our theorem, we need to recall two concepts. One
is  the Yamabe constant and the other is the best constant in the
Sobolev imbedding.

 Given a Riemannian metric $g(t)$ the Yamabe constant is
\[
Y(t) \equiv inf \frac{\int [|\nabla \phi|^2 + \frac{n-2}{4(n-1)} R
\phi^2] d\mu(x, t)}{\big{(} \int \phi^{2n/(n-2)} d\mu(x, t)
\big{)}^{(n-2)/n}}.
\]
The other is
 a Sobolev imbedding theorem due to E. Hebey
\cite{Heb:1} which is a refined form (on the controlling
constants) of the result by T. Aubin \cite{Au:1}:
\medskip

{\it {\bf Theorem S.} \ Let ${\bf M}$ be a complete (compact or
noncompact) Riemannian n-manifold. Suppose the Ricci curvature is
bounded below by $k$ and the injectivity radius is bounded below by
$i>0$. For any $\epsilon>0$, there exists $B(g)=B(\epsilon, n, k,
i)$  such that for any $\phi \in W^{1, 2}({\bf M})$,
\[
\bigg{(} \int_{\bf M}  |u|^{2n/(n-2)} d\mu(g) \bigg{)}^{(n-2)/n}
\le (K(n)^2 + \epsilon ) \int_{\bf M}  |\nabla u|^2 d\mu(g) + B(g)
\int_{\bf M}  u^2 d\mu(g).
\]Here $K(n)$ is the best constant in the Sobolev imbedding in ${\bf
R}^n$. }
\medskip

 We also need to mention the result by Hebey and Vaugon [HV]
where Theorem S is proven with $\epsilon=0$. However, then the
constant $B$ may depend on the derivative of the curvature tensor
which is harder to control. For our purpose, it suffices to fix the
$\epsilon$ as any positive constant, say $1$.

The following theorem is the main result of the section. It
contains three statements. The first one is an upper bound
controlled by the Yamabe constants, the second is an upper bound
controlled by the constant $B(g)$ in the Sobolev imbedding Theorem
S. They may seem technical at the first glance. However, the third
statement of the theorem provides a clarification. It shows that
these upper-bounds are the proper extension of on-diagonal upper
bound for the heat kernel in the fixed metric case. Recall that
for a compact Riemannian manifold ${\bf M}$ without boundary, the
heat kernel $G$ satisfies the following on-diagonal upper bound:
\[
G(x, t; y, s) \le c_1 \max \{ \frac{1}{(t-s)^{n/2}}, 1 \}
\]for some constant $c_1, c_2>0$ and for all $t>s$ and $x, y \in
{\bf M}$.

Here are some additional notations for the theorem. We will write
\[
R^- = - \min \{ R(x, t), 0 \}
\]where $R(x, t)$ is the scalar curvature under the metric $g(t)$.
When the scalar curvature changes sign, the theorem will also
involve the expression
\[
\frac{1}{ (\max R^-(\cdot, t))^{-1}+ (t-s)}.
\]This quantity is regarded as $0$ when $R(\cdot, t) \ge 0$.
\medskip

{\it Remark 4.1.} In statements (b) and (c) of Theorem 4.1 below,
the controlling constants depend only on the dimension, the lower
bound of Ricci curvature and the lower bound of injectivity radii.
By the result of Cheeger [Che], if one assumes that the sectional
curvatures are bounded between two constants and the volume of
geodesic balls of radius $1$ is bounded below by a positive
constant, then the injectivity radii are bounded from below by a
positive constant.  Therefore, the controlling constants in (b)
and (c) depend only on the bound of sectional curvature, the lower
bound of volume of balls of radius $1$ and dimension. The same can
be said for Theorem 5.1 below. The upshot is that the length of
time and the incompatibility of metric at different time do not
destroy the bound.

\medskip
\begin{theorem}
Suppose equation (4.0) has a smooth solution in the time interval
$[s, t]$ and let $G$ be the fundamental solution of (4.0). Then the
following statements hold.

(a). Suppose the Yamabe invariant $Y(g(\tau))>0$ for $\tau \in [s,
t]$, then
\[
\al
 G(x,& t; y, s) \\
 & \le  \frac{c_n}{\bigg{(} \int^{(t+s)/2}_s
 e^{2 c_n a(\tau)/n} Y(\tau) d\tau  \int^t_{(t+s)/2} e^{-2 c_n a(\tau)/n} \big{[} 1 +
(t-\tau)  \max R^-(\cdot, t) \big{]}^{-4/n} \ Y(\tau) d\tau
\bigg{)}^{n/4}}. \eal
\]Here $a(\tau) = \int^\tau_s \frac{1}{(\max R^-(\cdot, t))^{-1} + (t-l)}
dl.$

(b). Let $B(g(\tau))$ be the best constant in the Sobolev
imbedding Theorem S. Then
\[
G(x, t; y, s) \le \frac{c_n}{\bigg{(} \int^{(t+s)/2}_s e^{2
H(\tau) c_n /n} d\tau \ \int^t_{(t+s)/2} [ 1 + (t-\tau) \max
R^-(\cdot, t)]^{-4/n} e^{-2 H(\tau) c_n /n} d\tau \bigg{)}^{n/4}}
\]with
\[
H(\tau) =  \int^{\tau}_s [ B(g(l)) + \frac{1}{(\max{R^-(\cdot,
t)})^{-1} +  (t-l)}] dl.
\]

(c). In the special case that $R(\cdot, t) \ge 0$ and
$Ric(g(\tau)) \ge k$ and the injectivity radius is bounded below
by $i>0$, for all $\tau \in [s, t]$, then
\[
G(x, t; y, s) \le C(n, B) \max \{\frac{1}{(t-s)^{n/2}}, 1\}.
\]Here $B$ only depends on $n$, $k$ and $i$.  Moreover
\[
\frac{|\nabla_y G(x, t; y, s)|^2}{G(x, t; y, s)^2} \le C(n, B)
\frac{1}{(t-s)} \log \frac{\max \{\frac{1}{(t-s)^{n/2}}, 1\}}{G(x,
t; y, s)}.
\]
\end{theorem}
\medskip

{\it Remark 4.2.} Recently, in a paper [CL], Chang and Lu, proved
a derivative estimate for the Yamabe constant under the Ricci
flow. It can be coupled with this theorem to obtain better upper
bound on $G$.

Professor Lei Ni also informs us that he also knows a result on
upper bound in the case of certain Sobolev inequality.

\medskip

{\bf Proof of part (a).}  Without loss of generality, we take $s=0$
here and later.

Let $G$ be the fundamental solution to (4.0). By the reproducing
property
\[
G(x, t; y, 0) = \int G(x, t; z, t/2) G(z, t/2; y, 0) d\mu(z, t/2),
\]there holds
\[
G(x, t; y, 0) \le \bigg{[} \int G^2(x, t; z, t/2) d\mu(z, t/2)
\bigg{]}^{1/2} \ \bigg{[} \int G^2(z, t/2; y, 0) d\mu(z, t/2)
\bigg{]}^{1/2}. \leqno(4.1)
\]Therefore an upper bound follows from pointwise estimate on the
two quantities
\[
p(t) = \int G^2(x, t; y, s) d\mu(x, t), \leqno(4.2)
\]
\[
q(s) = \int G^2(x, t; y, s) d\mu(y, s). \leqno(4.3)
\]Let us estimate $p(t)$ in (4.2) first.

It is clear that
\[
\frac{d}{dt} p(t) = 2 \int G [ \Delta G - R G ] d\mu(x, t) + \int
G^2 R d\mu(x, t).
\]Here and later we omit the arguments on $G$ and differential
operators when no confusions appear. Therefore
\[
\al
 \frac{d}{dt} p(t) &\le - \int [ |\nabla G |^2 +  R G^2 ]
d\mu(x, t)\\
&= - \int [ |\nabla G |^2 +  \frac{n-2}{4(n-1)} R G^2 ] d\mu(x, t)
-\frac{3n-2}{4(n-1)} \int R G^2 d\mu(x, t). \eal
\leqno(4.4)
\]Let $Y(t)$ be the Yamabe constant with respect to $g(t)$, i.e.
\[
Y(t) = inf \frac{\int [|\nabla \phi|^2 + \frac{n-2}{4(n-1)} R
\phi^2] d\mu(x, t)}{\big{(} \int \phi^{2n/(n-2)} d\mu(x, t)
\big{)}^{(n-2)/n}}.
\]By H\"older's inequality
\[
\int G^2 d\mu(x, t) \le \big{[} \int G^{2n/(n-2)} d\mu(x, t)
\big{]}^{(n-2)/(n+2)} \ \big{[} \int G d\mu(x, t) \big{]}^{4/(n+2)},
\]we arrive at the 'conformal' Nash inequality
\[
\al & \int G^2 d\mu(x, t) \\
&\le c_n Y(t)^{-n/(n+2)} \ \big{[} \int [|\nabla \phi|^2 +
\frac{n-2}{4(n-1)} R \phi^2] d\mu(x, t) \big{]}^{n/(n+2)} \ \big{[}
\int G d\mu(x, t) \big{]}^{4/(n+2)}. \eal \leqno(4.5)
\]It is easy to check that
\[
\int G(x, t; y, s) d \mu(x, t) = 1.
\]Therefore (4.5) becomes
\[
\int G^2 d\mu(x, t) \\
c_n \le c_n Y(t)^{-n/(n+2)} \ \big{[} \int [|\nabla \phi|^2 +
\frac{n-2}{4(n-1)} R \phi^2] d\mu(x, t) \big{]}^{n/(n+2)}.
\]Substituting this to (4.4), we deduce
\[
p'(t) \le - c_n p(t)^{(n+2)/n} Y(t)- \frac{3n-2}{4(n-1)} \int R
G^2 d\mu(x, t).
\]It well-known (see [CK] e.g.) that the scalar curvature $R$ satisfies
the inequality
\[
\frac{d R}{dt} + \Delta R + c_n R^2 \le 0.
\]This implies
\[
R(y, \tau) \ge - \frac{1}{(\max R^-(\cdot, t))^{-1} + c_n
(t-\tau)}, \qquad \tau < t.
\]Here and later, if $R (\cdot, t) \ge 0$, then the above fraction is
 regarded
as zero.
Hence, for $\tau \in (s, t)$,
\[
p(\tau)' \le - c_n p(\tau)^{(n+2)/n} Y(\tau) + \frac{c_n}{(\max
R^-(\cdot, \tau))^{-1} + (t-\tau)} p(\tau).
\]Let
\[
a(\tau) = \int^\tau_s \frac{1}{(\max R^-(\cdot, l))^{-1} + (t-l)}
dl.
\]Then the above ordinary differential inequality becomes
\[
(e^{-c_n a(\tau)} p(\tau))' \le c_n \big{(} e^{-c_n a(\tau)}
p(\tau) \big{)}^{(n+2)/n} e^{2 c_n a(\tau)/n}  Y(\tau).
\]Integrating from $s$ to $t$, we deduce
\[
e^{-c_n a(t)} p(t) \le \frac{c_n}{\big{[} \int^t_s e^{2 c_n
a(\tau)/n} Y(\tau) d\tau \big{]}^{n/2}}
\]
This immediately shows that
\[
\int G^2(x, t; y, s) d\mu(x, t) = p(t) \le \frac{c_n e^{c_n
a(t)}}{\big{[} \int^t_s e^{2 c_n a(\tau)/n} Y(\tau) d\tau
\big{]}^{n/2}}. \leqno(4.6)
\]
\medskip

Next we estimate $q(s)$ in (4.3). Due to the asymmetry of the
equation, the computation is different.  Notice that the second
entries of $G$ satisfies the backward heat equation. i.e.
\[
\Delta_y G(x, t; y, s) + \partial_s G(x, t; y, s) =0.
\]This gives
\[
q'(s) = - 2 \int G \Delta G d\mu(y, s) + \int R G^2 d\mu(y, s).
\]Hence
\[
q'(s) \ge   \int [ |\nabla G|^2 + R G^2] d\mu(y, s). \leqno(4.7)
\]By the same argument as before we arrive at the Nash inequality
\[
\al & \int G^2 d\mu(y, s) \\
&\le c_n Y(s)^{-n/(n+2)} \ \big{[} \int [|\nabla \phi|^2 +
\frac{n-2}{4(n-1)} R \phi^2] d\mu(y, s) \big{]}^{n/(n+2)} \ \big{[}
\int G d\mu(y, s) \big{]}^{4/(n+2)}. \eal \leqno(4.8)
\]

This time we have to compute the quantity
\[
I(s) \equiv  \int G(x, t; y, s) d\mu(y, s).
\]It is clear that
\[
I'(s) = \int G(x, t; y, s) R(y, s) d\mu(y, s). \leqno(4.9)
\]Recall that
\[
R(y, \tau) \ge - \frac{1}{(\max R^-(\cdot, t))^{-1} + c_n (t-\tau)},
\qquad \tau < t.
\]Combining this with (4.9) we deduce
\[
I'(\tau) \ge  - \frac{1}{(\max R^-(\cdot, t))^{-1} + c_n (t-\tau)}
I(\tau).
\]Integrating from $s$ to $t$ and noting that $I(t) = 1$, we
obtain
\[
I(s) \le 1 + c_n (t-s) \max R^-(\cdot, t). \leqno(4.10)
\]Substituting (4.10) to (4.8), we deduce
\[
\al
 \int &[|\nabla G|^2 + R G^2] d\mu(y, s) \\
 &=\int [|\nabla G|^2 + \frac{n-2}{4(n-1)} R G^2] d\mu(y, s)
 + \frac{3n-2}{4(n-1)} \int  R G^2 d\mu(y, s)\\
 &\ge [ q(s)
]^{(n+2)/n} \big{[} 1 + c_n (t-s) \max R^-(\cdot, t)
\big{]}^{-4/n} \ Y(s)- \frac{c_n q(s)}{(\max R^-(\cdot, t))^{-1} +
(t-s)}. \eal \leqno(4.11)
\]Here, again we used the lower bound on the scalar curvature,
given just below (4.9).

This and (4.7) together imply that, for $\tau \in (s, t)$,
\[
q'(\tau) \ge c_n [ q(\tau) ]^{(n+2)/n} \big{[} 1 +  (t-\tau) \max
R^-(\cdot, t) \big{]}^{-4/n} \ Y(\tau) - \frac{c_n }{(\max
R^-(\cdot, t))^{-1} + (t-\tau)} q(\tau).
\]

Let again
\[
a(\tau) = \int^\tau_s \frac{1}{(\max R^-(\cdot, t))^{-1} + (t-l)}
dl.
\]Then
\[
\bigg{[} q(\tau) e^{c_n a(\tau)} \bigg{]}' \ge c_n \bigg{[}
q(\tau) e^{c_n a(\tau)} \bigg{]}^{(n+2)/n} \ e^{-2c_n a(\tau)/n}
\big{[} 1 +  (t-\tau) \max R^-(\cdot, t) \big{]}^{-4/n} \ Y(\tau).
\]
Integrating from $s$ to $t$, we obtain
\[
\al
 q(s)& = \int G^2(x, t; y, s) d\mu(y, s)\\
 & \le \frac{c_n e^{-c_n
a(s)} }{\bigg{(} \int^t_s e^{-2 c_n a(\tau)/n} \big{[} 1 +
(t-\tau) \max R^-(\cdot, t) \big{]}^{-4/n} \ Y(\tau) d\tau
\bigg{)}^{n/2}}. \eal
 \leqno(4.12)
\]

Now (4.6) an (4.12) respectively imply that
\[
\int G^2(z, t/2; y, 0) d\mu(z, t/2)  \le \frac{c_n e^{c_n
a(t/2)}}{\bigg{(} \int^{t/2}_0 e^{2 c_n a(\tau)/n}
 Y(\tau) d\tau \bigg{)}^{n/2}}, \leqno(4.13)
\]
\[
\al
 \int& G^2(x, t; z, t/2) d\mu(z, t/2)\\
 &  \le \frac{c_n e^{-c_n
a(t/2)} }{\bigg{(} \int^t_{t/2} e^{-2 c_n a(\tau)/n} \big{[} 1 +
(t-\tau) \max R^-(\cdot, t) \big{]}^{-4/n} \ Y(\tau) d\tau
\bigg{)}^{n/2}}, \eal
 \leqno(4.14)
\]By (4.1), (4.13) and (4.14), we arrive at the following upper bound
\[
\al
 G(x,& t; y, 0) \\
 & \le  \frac{c_n}{\bigg{(} \int^{t/2}_0
 e^{2 c_n a(\tau)/n} Y(\tau) d\tau  \int^t_{t/2} e^{-2 c_n a(\tau)/n} \big{[} 1 +
(t-\tau)  \max R^-(\cdot, t) \big{]}^{-4/n} \ Y(\tau) d\tau
\bigg{)}^{n/4}}. \eal
\]This proves part (a).

\medskip

{\bf Proof of Part (b).}

We generally follow the previous arguments  between (4.1) and
(4.14) to derive an upper bound. The difference is that we will
use the Sobolev inequality (Theorem S) instead of the Yamabe
constant.

As before, by H\"older's inequality and Theorem S, we arrive at
the Nash type inequality
\[
\al & \int G^2 d\mu(x, t) \\
&\le  \ \big{[} \int [c_n |\nabla G|^2 d\mu(x, t) + B(g(t))
\int_{\bf M} G^2 d\mu(x, t) \big{]}^{n/(n+2)} \ \big{[} \int G
d\mu(x, t) \big{]}^{4/(n+2)}. \eal \leqno(4.15)
\]Here and later $c_n$ is a dimensional constant that may change
from line and to line. Since, again,
\[
\int G d \mu(x, t) = 1,
\]we have
\[
\int |\nabla G|^2 d\mu(x, t) \ge c_n \bigg{[} \int_{\bf M} G^2
d\mu(x, t)  \bigg{]}^{(n+2)/n} - c_n B(g(t)) \int G^2 d\mu(x, t).
\leqno(4.16)
\]Combing (4.16) with (4.4) under again the notation (4.2), we
obtain
\[
p'(t) \le -c_n p(t)^{(n+2)/n} + c_n B(g(t)) p(t) - \int_{\bf M} R
G^2 d\mu(x, t).
\]Fixing $s$ and $t$, for any $\tau \in (s, t)$, we still have the lower
bound for the scalar curvature (just after (4.9))
\[
R(\cdot, \tau) \ge - \frac{1}{(\max R^-(\cdot, t))^{-1} + c_n
(t-\tau)}, \qquad \tau < t.
\]Therefore
\[
p'(\tau) \le -c_n p(\tau)^{(n+2)/n} + c_n h(\tau) p(\tau),
\leqno(4.17)
\]with
\[
h(\tau) =  B(g(\tau)) + \frac{1}{(\max{R^-(\cdot, t)})^{-1} +
(t-\tau)}.
\]

Let $ H(\tau)$ be the anti-derivative of $h(\tau)$ such that $H(s)
= 0$. Then
\[
\bigg{(} e^{-c_n H(\tau)} p(\tau) \bigg{)}' \le - c_n \bigg{(}
e^{- c_n H(\tau)} p(\tau) \bigg{)}^{(n+2)/n} \ e^{2 c_n
H(\tau)/n}.
\]Integrating from $s$ to $t$, we arrive at
\[
p(t) \le  \frac{c_n e^{c_n H(t)}}{\bigg{(} \int^t_s e^{2 H(\tau)
c_n /n} d\tau \bigg{)}^{n/2}}. \leqno(4.18)
\]

Our next task is to bound
\[
q(s) = \int G^2(x, t; y, s) d\mu(y, s).
\]Clearly the counter-parts of (4.7), (4.10) and (4.15) still
hold. i.e.
\[
q'(s) \ge \int (|\nabla G|^2 + R G^2) d\mu(y, s).
\]
\[
I(s) = \int G d\mu(y, s) \le 1 + c_n (t-s) \max R^-(\cdot, t).
\]
\[
\al & \int G^2 d\mu(y, s) \\
&\le  \ \big{[} \int [c_n |\nabla G|^2 d\mu(y, s) + B(g(s)) \int
G^2 d\mu(x, t) \big{]}^{n/(n+2)} \ \big{[} \int G d\mu(y, s)
\big{]}^{4/(n+2)}. \eal
\]Also
\[
R(y, s) \ge - \frac{1}{[\max R^-(\cdot, t)]^{-1} + c_n (t-s)}.
\]These four inequalities imply that
\[
q'(s) \ge c_n q(s)^{(n+2)/n} [ 1 + (t-s) \max R^-(\cdot,
t)]^{-4/n} - c_n h(s) q(s).
\]Here $h(s)$ is given by the expression just below (4.17) with
$\tau$ replaced by $s$. Now, for fixed $s$ and $t$ and any $\tau
\in (s, t)$, the above differential inequality on $q'(s)$ is still
valid for $q'(\tau)$ when $s$ is replaced by $\tau$. Let $H(\tau)$
be the antiderivative of $h(\tau)$ with $H(s) = 0$. Then it is
clear that
\[
\bigg{(} e^{c_n H(\tau)} q(\tau) \bigg{)}' \ge c_n \bigg{(} e^{c_n
H(\tau)} q(\tau) \bigg{)}^{(n+2)/n} [ 1 + (t-\tau) \max R^-(\cdot,
t)]^{-4/n} e^{-2 H(\tau) c_n /n}.
\]Integrating from $s$ to $t$, we arrive at
\[
q(s) \le  \frac{c_n e^{-c_n H(s)}}{\bigg{(} \int^t_s [ 1 +
(t-\tau) \max R^-(\cdot, t)]^{-4/n} e^{-2 H(\tau)c_n/n} d\tau
\bigg{)}^{n/2}}. \leqno(4.19)
\]

By (4.18), we have
\[
p(t/2)= \int G^2(z, t/2; y, 0) d\mu(z, t/2) \le \frac{c_n e^{c_n
H(t/2)}}{\bigg{(} \int^{t/2}_0 e^{2 H(\tau) c_n /n} d\tau
\bigg{)}^{n/2}}. \leqno(4.20)
\]

Also, (4.19) shows
\[
q(t/2)=\int G^2(x, t; z, t/2) d\mu(z, t/2) \le \frac{c_n
e^{-H(t/2)}}{\bigg{(} \int^t_{t/2} [ 1 + (t-\tau) \max R^-(\cdot,
t)]^{-4/n} e^{-2 H(\tau)/n} d\tau \bigg{)}^{n/2}}. \leqno(4.21)
\]Here
\[
H(t/2) = \int^{t/2}_0 [ B(g(\tau)) + \frac{1}{(\max{R^-(\cdot,
t)})^{-1} +  (t-\tau)}] d\tau.
\]

Multiplying (4.20) and (4.21), and using (4.1), we have proven the
on-diagonal upper bound
\[
G(x, t; y, 0)^2 \le \frac{c_n}{\bigg{(} \int^{t/2}_0 e^{2 H(\tau)
c_n /n} d\tau \bigg{)}^{n/2} \ \bigg{(} \int^t_{t/2} [ 1 +
(t-\tau) \max R^-(\cdot, t)]^{-4/n} e^{-2 H(\tau) c_n/n} d\tau
\bigg{)}^{n/2}}
\]This gives part (b).
\medskip

{\bf Proof of part (c).}

In the special case that $R(\cdot, t) \ge 0$ and $Ricc(g(\tau))
\ge k$ uniformly and the injectivity radius is uniformly bounded
below by $i$, then
\[
H(\tau) = \int^{\tau}_0 [ B(g(l)) + \frac{1}{(\max{R^-(\cdot,
t)})^{-1} +  (t-l)}] dl =  c_n B(n, k, i) \tau
\]and
\[
 \frac{1}{(\max{R^-(\cdot,
t)})^{-1} +  (t-\tau)}=0.
\]Hence the above immediately shows
\[
G(x, t; y, 0) \le C(c_n, B) \max \{\frac{1}{t^{n/2}}, 1\}.
\]

The gradient estimate follows from Hamilton's argument in Theorem
1.1 [H], which can be easily generalized to the present case.  For
$z \in {\bf M}$ and $\tau \in [0, t/2]$, let
\[
v(z, \tau) = G(x, t; z, \tau).
\]Then $v$ is a solution to the backward heat equation $\Delta v +
v_{\tau} = 0$. By direct computation
\[
\al (\Delta +\partial_\tau)&(\frac{|\nabla v|^2}{v}) = \frac{2}{v}
\bigg{|}\partial_i \partial_j v - \frac{\partial_i v \partial_j
v}{v} \bigg{|}^2\\
& \ge 0.
 \eal
\]

Let $A$ be the maximum of $v$ in the time interval $[0, t/2]$. By
the above estimate
\[
A \le C(c_n, B) \max \{\frac{1}{t^{n/2}}, 1\}.
\]Direct computation shows
\[
\begin{aligned}
 &\Delta (v \log \frac{A}{v}) + \partial_\tau (v \log \frac{A
}{v}) \\
&=(\Delta v + \partial_\tau v) \log A + (\Delta + \partial_\tau) (v
\log\frac{A}{v})\\
&=-\frac{|\nabla v|^2}{v}.
 \end{aligned}
\]Let $\phi = (t/2)-\tau$, then it is clear that, for
\[
h = \phi \frac{|\nabla v|^2}{v} - v \log \frac{A }{v},
\]there holds
\[
\Delta h + \partial_\tau h \ge 0.
\]By the maximum principle, applied backward in time, we have
\[
\frac{|\nabla_y v|^2}{v^2} \le C \frac{1}{\tau} \log \frac{\max
\{\frac{1}{t^{n/2}}, 1\}}{v}
\]for $\tau \in [0, t/4]$. \qed

\section{The case of nonnegative Ricci curvature}

In this section, we specialize to the case of nonnegative Ricci
curvature. We establish certain Gaussian type upper bound for the
fundamental solution of (4.0). We will begin with the traditional
method of establishing a mean value inequality via Moser's iteration
and a weighted estimate in the spirit of Davies [Da]. However, there
is some difficulty in applying this method directly due to the lack
of control of the time derivative of the distance function. The new
idea to overcome this difficulty is to use the interpolation result
of Theorem 3.2 and the bound in Theorem 4.1 (c).

The following are some additional notations for this section. We
will use $B(x, r; t)$ to denote the geodesic ball centered at $x$
with radius $r$ under the metric $g(t)$; $|B(x, r; t)|_s$ to denote
the volume of $B(x, r; t)$ under the metric $g(s)$.

The main result of this section is Theorem 5.1 below. Note that the
theorem is qualitatively sharp in general since it matched the
well-known Gaussian upper bound for the fixed metric case. Also
there is no assumption on the comparability of metrics at different
times.  In this theorem, we assume the manifold is compact. This
accounts for the extra $1$ on the Gaussian upper bound. Even in the
case of fixed metric, the heat kernel converges to a positive
constant for large time. The theorem still holds for certain
noncompact manifolds under suitable assumptions. In this case the
extra $1$ in the upper bound should be replaced by $0$.

\medskip
{\it Remark 5.1.} As mentioned in section 4 (Remark 4.1), the
controlling constants in the theorem below can be made to depend
only on the bound of sectional curvature, the lower bound of
volume of balls of radius $1$ and the dimension.

In the case $Ricci \ge - k$ with $k>0$, then certain integral
Gaussian bound similar to the one below (5.13) can still be proven
by the same method. However, so far we are not able to derive a
pointwise Gaussian upper bound without an exponentially growing term
$e^{k t}$. This is due to a lack of an efficient mean value
inequality for the second entries of the fundamental solution, which
satisfies (3.2) after a time reversal.

\medskip
\begin{theorem}
Assume that equation (4.0) has a smooth solution in the time
interval $[0, T]$ and let $G$ be the fundamental solution of
(4.0). Suppose that $Ricci \ge 0$ and that the injectivity radius
is bounded from below by a positive constant $i$ throughout. Then
the following statement holds.

 For any $s, t \in (0, T)$ and $x, y \in {\bf M}$,  there exist a dimensional constant
$c_n$, a dimension less constant $c$ and a constant $A$ depending
only on $i$ such that
\[
G(x, t; y, s) \le c_n A  \bigg{(} 1+  \frac{1}{(t-s)^{n/2}} +
 \frac{1}{|B(x, \sqrt{t-s},
t)|_s}  \bigg{)}   e^{-c d(x, y, s)^2/(t-s)}.
\]
\end{theorem}

{\bf Proof}

It is obvious that we only have to deal with the case that $B(x,
2\sqrt{t-s}, s)$ is a proper sub-domain of ${\bf M}.$  Otherwise,
$\sqrt{t-s}/2 \ge  d(x, y, s)$ for any $x, y \in {\bf M}$. So the
exponential term is mute and the result is already proven by
Theorem 4.1 (c).

First we use Moser's iteration to prove a mean value inequality. The
only new factor is a cancelation effect induced by the backward
Ricci flow. So we will be brief in the presentation at this part of
the proof.

Let $u$ be a positive solution to (4.0) in the region
\[
Q_{\sigma r}(x, t) \equiv \{ (y, s) \ | \ z \in {\bf M}, t-(\sigma
r)^2 \le s \le t, \ d(y, x, s) \le \sigma r \}.
\]Here $r>0, 2 \ge \sigma \ge 1$.
Given any $p \ge 1$, it is clear that
\[
\Delta u^p - p R u^p- \partial_t u^p \ge 0. \leqno(5.1)
\]

Let $\phi: [0, \infty) \to [0, 1]$ be a smooth function such that
$|\phi'| \le 2/((\sigma-1) r)$, $\phi' \le 0$, $\phi \ge 0$,
$\phi(\rho) =1$ when $0 \le \rho \le r$, $\phi(\rho)=0$ when $\rho
\ge \sigma r$. Let $\eta: [0, \infty) \to [0, 1]$ be a smooth
function such that $|\eta'| \le 2/((\sigma-1) r)^2$, $\eta' \ge
0$, $\eta \ge 0$, $\phi(s) =1$ when $t-r^2 \le s \le t$,
$\phi(s)=0$ when $s \le t- (\sigma r)^2$.

Writing $w= u^p$ and using $w \psi^2$ as a test function on (5.2),
we deduce
\[
\int \nabla (w \psi^2) \nabla w d\mu(y, s)ds + p \int R w^2 \psi^2
d\mu(y, s)ds \le - \int (\partial_s w) w \psi^2  d\mu(y, s)ds.
\leqno(5.3)
\]By direct calculation
\[
\int \nabla (w \psi^2) \nabla w d\mu(y, s)ds = \int  | \nabla (w
\psi)|^2 d\mu(y, s)ds - \int  |\nabla \psi|^2 w^2 d\mu(y, s)ds.
\leqno(5.4)
\]Next we estimate the righthand side of (5.3). Here we will use the
backward Ricci flow.
\[
- \int (\partial_s w) w \psi^2  d\mu(y, s)ds = \int w^2 \psi
\partial_s \psi d\mu(y, s)ds + \frac 1 2 \int (w \psi)^2 R d\mu(y, s)ds
- \frac{1}{2} \int (w \psi)^2 d\mu(y, t).
\]Observe that
\[
\partial_s \psi = \eta(s)   \phi'(d(y, x, s)) \partial_s d(y, x, s)
+ \phi(d(y, x, s)) \eta'(s) \le \phi(d(y, x, s)) \eta'(s).
\]This is so because $\phi' \le 0$ and $\partial_s d(y, x, s) \ge 0$
under the backward Ricci flow with nonnegative Ricci curvature.
Hence
\[
\al
 - \int &(\partial_s w) w \psi^2  d\mu(y, s)ds \\
 &\le \int w^2 \psi \phi(d(y, x, s))  \eta'(s)  d\mu(y, s)ds +
 \frac 1 2 \int (w \psi)^2 R d\mu(y, s)ds
- \frac{1}{2} \int (w \psi)^2 d\mu(y, t). \eal \leqno(5.5)
\]Combing (5.3) to (5.5), we obtain, in view of $p \ge 1$ and $R \ge
0$,
\[
\int  | \nabla (w \psi)|^2 d\mu(y, s)ds + \frac{1}{2} \int (w
\psi)^2 d\mu(y, t) \le \frac{c}{(\sigma-1)^2 r^2} \int_{Q_{\sigma r
(x, t)}} w^2 d\mu(y, s)ds. \leqno(5.6)
\]By H\"older's inequality
\[
\int (\psi w)^{2(1+(2/n)} d\mu(y, s) \le \bigg{(}   \int (\psi
w)^{2n/(n-2))} d\mu(y, s) \bigg{)}^{(n-2)/n} \bigg{(} \int (\psi
w)^2 d\mu(y, s) \bigg{)}^{2/n}. \leqno(5.7)
\]
\medskip

Let us assume that  $B(x, \sigma r, s)$ is a proper sub-domain of
${\bf M}.$ In this case, for manifolds with nonnegative Ricci
curvature, it is well-known that the following Sobolev imbedding
holds (see [Sa] e.g.)
\[
\bigg{(}  \int (\psi w)^{2n/(n-2)} d\mu(y, s) \bigg{)}^{(n-2)/n}
\le \frac{c_n \sigma^2 r^2}{|B(x, \sigma r, s)|^{2/n}_s}  \int [|
\nabla (\psi w) |^2  + r^{-2} (\psi w)^2] d\mu(y, s).
\]For $s \in [t-(\sigma r)^2, t]$, by the assumption that the Ricci
curvature is nonnegative, it holds
\[
B(x, \sigma r, s) \supset B(x,  \sigma r, t); \qquad |B(x, \sigma r,
s)|_s \ge |B(x, \sigma r, t)|_{t-(\sigma r)^2}.
\]Therefore we have
\[
\bigg{(}  \int (\psi w)^{2n/(n-2)} d\mu(y, s) \bigg{)}^{(n-2)/n}
\le \frac{c_n \sigma^2 r^2}{|B(x, \sigma r, t)|^{2/n}_{t-(\sigma
r)^2}} \int [| \nabla (\psi w) |^2  + r^{-2} (\psi w)^2] d\mu(y,
s). \leqno(5.8)
\]For $s \in [t-(\sigma r)^2, t]$. Substituting (5.7) and (5.8) to
(5.6), we arrive at the estimate
\[
\int_{Q_r(x, t)} w^{2 \theta} d\mu(y, s) ds \le c_n \frac{r^2}{|B(x,
\sigma r, t)|^{2/n}_{t-(\sigma r)^2}} \bigg{(} \frac{1}{(\sigma-1)^2
r^2} \int_{Q_{\sigma r}(x, t)} w^2 d\mu(y, s) ds \bigg{)}^{\theta},
\]with $\theta = 1+(2/n)$. Now we apply the above inequality with
the parameters $\sigma_0=1, \sigma_i=1- \Sigma^i_{j=1} 2^{-j-1}$ and
$p=\theta^i$. This shows a $L^2$ mean value inequality
\[
\sup_{Q_{r/2}(x, t)} u^2 \le \frac{c_n}{r^2 |B(x, r, t)|_{t-r^2}}
\int_{Q_r(x, t)} u^2 d\mu(y, s)ds. \leqno(5.9)
\]From here, by a generic trick of Li and Schoen [LS], applicable here since it
 uses
only the doubling property of the metric balls, we arrive at the
$L^1$ mean value inequality
\[
\sup_{Q_{r/2}(x, t)} u \le \frac{c_n}{r^2 |B(x, r, t)|_{t-r^2}}
\int_{Q_r(x, t)} u d\mu(z, \tau)d\tau. \leqno(5.10)
\]

Fixing $y \in {\bf M} $ and $s<t$, we apply (5.10) on $u = G(\cdot,
\cdot; y, s)$ with $r=\sqrt{t-s}/2$. Note that $\int_{\bf M} u(z,
\tau) d\mu(z, \tau) =1$. The doubling property of the geodesic balls
show that
\[
G(x, t; y, s) \le \frac{c_n}{|B(x, \sqrt{t-s}, t)|_s}
\]when $|B(x, \sqrt{t-s}, s)|$ is a proper subdomain of ${\bf M}$.

\medskip

Without loss of generality, we take $s=0$.
 We begin by using a
modified version of the exponential weight method due to Davies
[Da]. Pick a point $x_0 \in {\bf M}$, a number $\lambda<0$ and a
function $f \in L^2({\bf M}, g(0))$. Consider the functions $F$ and
$u$ defined by
\[
F(x, t) \equiv  e^{\lambda d(x, x_0, t)} u(x, t) \equiv e^{\lambda
d(x, x_0, t)} \int G(x, t; y, 0) e^{-\lambda d(y, x_0, 0)} f(y)
d\mu(y, 0).
\]It is clear that $u$ is a solution of (4.0).  By direct
computation, we have
\[
\al \partial_t &\int F^2(x, t) d\mu(x, t) = \partial_t \int
e^{2\lambda d(x, x_0, t)} u^2(x, t)d\mu(x, t)\\
&= 2\lambda \int e^{2\lambda d(x, x_0, t)} \partial_t d(x, x_0, t)
u^2(x, t)d\mu(x, t) + \int e^{2\lambda d(x, x_0, t)} u^2(x, t) R(x,
t) d\mu(x, t)\\
&\qquad + 2 \int e^{2\lambda d(x, x_0, t)}[\Delta u - R(x, t) u(x,
t)] u(x, t)d\mu(x, t).
 \eal
\]By the assumption that $Ricci \ge 0$ and $\lambda<0$, the above
shows
\[
\partial_t  \int F^2(x, t) d\mu(x, t)
\le 2  \int e^{2\lambda d(x, x_0, t)} u \Delta u d\mu(x, t).
\]Using integration by parts, we turn the above inequality into
\[
\al
\partial_t  &\int F^2(x, t) d\mu(x, t)\\
&\le -4  \lambda \int e^{2\lambda d(x, x_0, t)} u \nabla d(x, x_0,
t) \nabla u
 d\mu(x, t) -2 \int e^{2\lambda d(x, x_0, t)} |\nabla u|^2 d\mu(x,
t). \eal
\]Observe also
\[
\al \int &|\nabla F(x, t)|^2 d\mu(x, t) = \int |\nabla (e^{\lambda
d(x, x_0, t)} u(x, t))|^2 d\mu(x, t)\\
&=\int e^{2\lambda d(x, x_0, t)} |\nabla u|^2 d\mu(x, t) + 2 \lambda
\int e^{2\lambda d(x, x_0, t)} u \nabla d(x, x_0, t) \nabla u
 d\mu(x, t)\\
 &\qquad + \lambda^2 \int e^{2\lambda d(x, x_0, t)} |\nabla d|^2 u^2
 d\mu(x, t).
 \eal
 \]Combining the last two expressions, we deduce
 \[
\partial_t  \int F^2(x, t) d\mu(x, t)
\le - 2 \int |\nabla F(x, t)|^2 d\mu(x, t) + \lambda^2 \int
e^{2\lambda d(x, x_0, t)} |\nabla d|^2 u^2
 d\mu(x, t).
 \]By the definition of $F$ and $u$, this shows
 \[
\partial_t  \int F^2(x, t) d\mu(x, t) \le \lambda^2 \int F(x, t)^2
 d\mu(x, t).
\]Upon integration, we derive the following $L^2$ estimate
\[
\int F^2(x, t) d\mu(x, t) \le e^{\lambda^2 t} \int F^2(x, 0) d\mu(x,
0) = e^{\lambda^2 t} \int f(x)^2 d\mu(x, 0). \leqno(5.11)
\]

Recall that $u$ is a solution to (4.0). Therefore, by the mean value
inequality (5.9), the following holds
\[
u(x, t)^2 \le \frac{c_n}{t |B(x, \sqrt{t/2}, t)|_{t/2}} \int^t_{t/2}
\int_{B(x, \sqrt{t/2}, \tau)} u^2(z, \tau) d\mu(z, \tau) d\tau.
\]i.e. By the definition  of $F$ and $u$, it follows that
\[
u(x, t)^2 \le \frac{c_n}{t |B(x, \sqrt{t/2}, t)|_{t/2}} \int^t_{t/2}
\int_{B(x, \sqrt{t/2}, \tau)} e^{-2 \lambda d(z, x_0, \tau)} F^2(z,
\tau) d\mu(z, \tau) d\tau.
\]In particular, this holds for $x=x_0$. In this case, for $z \in
B(x_0, \sqrt{t/2}, \tau)$, there holds $d(z, x_0, \tau) \le
\sqrt{t/2}.$ Therefore, by the assumption that $\lambda<0$,
\[
u(x_0, t)^2 \le \frac{c_n e^{-2 \lambda \sqrt{t/2}}}{t |B(x_0,
\sqrt{t/2}, t)|_{t/2}} \int^t_{t/2} \int_{B(x_0, \sqrt{t/2}, \tau)}
F^2(z, \tau) d\mu(z, \tau) d\tau.
\]This combined with (5.11) shows that
\[
u(x_0, t)^2 \le \frac{c_n e^{ \lambda^2 t- \lambda \sqrt{2t}}}{
|B(x_0, \sqrt{t/2}, t)|_{t/2}}  \int f(y)^2 d\mu(y, 0).
\]i.e.
\[
\bigg{(} \int G(x_0, t; z, 0) e^{-\lambda d(z, x_0, 0)} f(z) d\mu(z,
0) \bigg{)}^2 \le \frac{c_n e^{ \lambda^2 t- \lambda \sqrt{2t}}}{
|B(x_0, \sqrt{t/2}, t)|_{t/2}}  \int f(y)^2 d\mu(y, 0). \leqno(5.12)
\]Now, we fix $y_0$ such that $d(y_0, x_0, 0)^2 \ge 4 a^2 t$ with $a>1$ to be chosen
later. Then it is clear that, by $\lambda<0$ and the triangle
inequality,
\[
-\lambda d(z, x_0, 0) \ge - a \lambda d(x_0, y_0, 0)
\]when $d(z, y_0, 0) \le \sqrt{t}$.
In this case, (5.12) implies
\[
\bigg{(} \int_{B(y_0, \sqrt{t}, 0)}  G(x_0, t; z, 0) f(z) d\mu(z, 0)
\bigg{)}^2 \le \frac{c_n e^{2 a \lambda d(x_0, y_0, 0) + \lambda^2
t- \lambda \sqrt{2t})}}{ |B(x_0, \sqrt{t/2}, t)|_{t/2}}  \int f(y)^2
d\mu(y, 0). \leqno(5.13)
\]Now we take
\[
\lambda = - \frac{d(x_0, y_0, 0)}{bt}.
\]Take $b>0$ and $a>0$ sufficiently large.
Then (5.13) shows, for some $c>0$,
\[
\int_{B(y_0, \sqrt{t}, 0)}  G^2(x_0, t; z, 0) d\mu(z, 0) \le
\frac{c_n e^{-c d(x_0, y_0, 0)^2/t}}{ |B(x_0, \sqrt{t/2},
t)|_{t/2}}.
\]Hence, there exists $z_0 \in B(y_0, \sqrt{t}, 0)$ such that
\[
G^2(x_0, t; z_0, 0) \le \frac{c_n e^{-c d(x_0, y_0, 0)^2/t}}{
|B(x_0, \sqrt{t/2}, t)|_{t/2} ||B(x_0, \sqrt{t}, 0)|_0}.
\]By the doubling property of the geodesic balls, it implies
\[
G^2(x_0, t; z_0, 0) \le \frac{c_n e^{-c d(x_0, y_0, 0)^2/t}}{
|B(x_0, \sqrt{t}, t)|_{0} ||B(x_0, \sqrt{t}, 0)|_0}. \leqno(5.14)
\]

Finally, let us remind ourself that $G(x_0, t; \cdot, \cdot)$ is a
solution to the conjugate equation of (4.0). i.e.
 \[ \Delta_z G(x,
t; z; \tau) + \partial_\tau G(x, t; z, \tau) =0. \] Therefore
Theorem 3.2 can be applied to it after a reversal in time.
Consequently, for $\delta>0, C>0$,
\[
G(x_0, t; y_0, 0) \le C G^{1/(1+\delta)} (x_0, t, z_0, 0)
M^{\delta/(1+\delta)}, \leqno(5.15)
\]where $M =\sup_{M \times [0, t/2]} G(x_0, t, \cdot, \cdot)$.
By Theorem 4.1, part (c), there exists a constant $A>0$, depending
only on the lower bound of the injectivity radius such that
\[
M \le A \max \{\frac{1}{t^{n/2}}, 1 \}.
\]This, (5.14) and (5.15) show, with $\delta = 1$, that
\[
G(x_0, t; y_0, 0)^2 \le \max \{ \frac{c_n}{t^{n/2}}, 1 \} \frac{A
e^{-c d(x_0, y_0, 0)^2/t}}{ \sqrt{|B(x_0, \sqrt{t}, t)|_{0}
||B(x_0, \sqrt{t}, 0)|_0}}.
\]
By the assumption that the Ricci curvature is nonnegative, we have
\[
 |B(x_0, \sqrt{t}, t)|_0 \le |B(x_0, \sqrt{t}, 0)|_0.
\]Therefore
\[
G^2(x_0, t; y_0, 0) \le \max \{ \frac{c_n}{t^{n/2}}, 1 \} \frac{A
e^{-c d(x_0, y_0, 0)^2/t}}{|B(x_0, \sqrt{t}, t)|_0}.
\]Consequently
\[
G(x_0, t; y_0, 0) \le c_n A  \bigg{(} 1+  \frac{1}{t^{n/2}} +
 \frac{1}{|B(x_0, \sqrt{t},
t)|_0}  \bigg{)}   e^{-c d(x_0, y_0, 0)^2/t}.
\]Since $x_0$ and $y_0$ are arbitrary, the proof is done.

 \qed

\bigskip

{\bf Acknowledgement.} We thank Professors Bennet Chow, Peng Lu and
Lei Ni and also Shilong Kuang for very helpful communications.

\bigskip

\noindent e-mail:  qizhang@math.ucr.edu

\enddocument